\documentclass[12pt]{article}
\usepackage{amssymb,amsmath,amsfonts}
\begin{document}

\newtheorem{teo}{Theorem}
\newtheorem{cor}{Corollary}
\newtheorem{lm}{Lemma}
\newtheorem{rem}{Remark}
\renewcommand{\theteo}{\arabic{teo}.}
\renewcommand{\thecor}{\arabic{cor}.}
\renewcommand{\thelm}{\arabic{lm}.}
\renewcommand{\therem}{\arabic{rem}.}
\newcommand{\indlim}{\operatornamewithlimits{ind\ lim}}
\renewcommand{\refname}{\begin{center} \mdseries {\sf \large References}
\end{center}}
\begin{center}{\bf ON THE CAUCHY PROBLEM FOR DIFFERENTIAL EQUATIONS IN
A BANACH SPACE OVER THE FIELD OF $p$-ADIC NUMBERS. I.}
\footnote{
Supported by CRDF (Project UM 1-2421-KV-02)}
\end{center}
\vspace*{3mm}
\begin{center} {\textmd \sf MYROSLAV L. GORBACHUK and VALENTYNA I. GORBACHUK}
\end{center}
\vspace*{3mm}
\begin{quote} ABSTRACT. \ For the Cauchy problem for an operator differential
equation of the form $y'(z) = Ay(z)$, where $A$ is a closed linear operator
on a Banach space over the field of $p$-adic numbers, the criterion of
well-posedness in the class of locally analytic vector-functions is established.
It is shown how the Cauchy-Kovalevskaya theorem for $p$-adic partial
differential equations may be obtained as a particular case from this
criterion.
\end{quote}
\vspace*{3mm}

{\bf 1.} \ Let $\mathfrak B$ be a Banach space with norm $\|\cdot\|$ over
the completion $\Omega = {\Omega}_p$ of an algebraic closure of the field of
$p$-adic numbers [1 - 3] ($p$ is prime),
and let $A$ be a closed linear operator on $\mathfrak B$, that is,
the convergences ${\mathcal D}(A) \ni x_n \to x$ and $Ax_n \to y \ (n \to \infty)$
in $\mathfrak B$ imply the inclusion $x \in {\mathcal D}(A)$ and the equality
$Ax = y \ ({\mathcal D}(\cdot)$ is the domain of an operator).

For a number $\alpha > 0$, we put
$$
E_{\alpha}(A) = \bigl\{x \in
\bigcap\limits_{n \in {\mathbb N}_0 = \{0, 1, 2, \dots\}}{\mathcal D}(A^n) \Big|
\exists c = c(x) > 0 \quad \forall k \in {\mathbb N}_0 \quad
\|A^k x\| \le c{\alpha}^k\bigr\}.
$$
The linear space $E_{\alpha}(A)$ is a Banach space with respect to the norm
$$
\|x\|_{\alpha} = \sup\limits_{n \in {\mathbb N}_0}\frac{\|A^nx\|}{{\alpha}^n}.
$$
The set
$$
E(A) = \bigcup\limits_{\alpha > 0} E_{\alpha}(A)
$$
is endowed with the inductive limit topology of the Banach spaces
$E_{\alpha}(A)$:
$$
E(A) = \indlim \limits_{\alpha \to \infty} E_{\alpha}(A).
$$
By the closed graph theorem, $E(A)$ coincides with $\mathfrak B$ if and only if
${\mathcal D}(A) = \mathfrak B$. The elements of $E(A)$ are called entire vectors of
exponential type for the operator $A$. Define  the type $\sigma (x; A)$ of
a vector $x \in E(A)$ as
$$
\sigma (x; A) = \inf \{\alpha > 0: x \in E_{\alpha}(A)\} =
\varlimsup_{n \to \infty} \|A^nx\|^{\frac{1}{n}}.
$$
Thus, the equality $\sigma (x; A) = \sigma$ means that for an arbitrary
$\varepsilon > 0$, there exists a constant $c_{\varepsilon} =
c_{\varepsilon}(x) > 0$ such that
$$
\forall n \in {\mathbb N}_0 \quad \|A^nx\| \le c_{\varepsilon}(\sigma + \varepsilon)^n,
$$
and
$$
\lim\limits_{i \to \infty} \frac{\|A^{n_i}x\|}{(\sigma - \varepsilon)^{n_i}} =
\infty
$$
for some subsequence $n_i \to \infty$ when $i \to \infty$. In the case, where
the operator $A$ is bounded, the type of a vector $x \in {\mathfrak B} = E(A)$
does not exceed the norm of $A$:
$$
\forall x \in {\mathfrak B} \quad \sigma (x; A) \le \|A\|.
$$

{\bf 2.} \ The object of consideration now is a power series
$$
y(z) = \sum\limits_{n = 0}^{\infty} c_nz^n, \quad c_n \in {\mathfrak B},
\quad z \in \Omega.  \eqno (1)
$$
For such a series the convergence radius is determined by the formula
$$
r = r(y) = \frac{1}{\varlimsup_{n \to \infty}\sqrt[n]{\|c_n\|}}. \eqno (2)
$$
In the open disk $d(r^-; \Omega) = \{z \in \Omega: |z|_p < r\}$, the series (1)
defines a vector-function $y(z)$ with values in $\mathfrak B$ ($|\cdot|_p$ is
the $p$-adic valuation on $\Omega$).

Denote by ${\mathfrak A}_{loc}(\mathfrak B)$ the set of all vector-functions $y(z)$
which are represented by a series of the form (1) with $r(y) > 0$. It is obvious
that ${\mathfrak A}_{loc}(\mathfrak B)$ is a vector space over $\Omega$. We
call its elements locally analytic vector-functions. The convergence $y_n \to
y \ (n \to \infty)$ in ${\mathfrak A}_{loc}(\mathfrak B)$ means that there exists
a number $\delta > 0$ such that $r(y_n) \ge \delta$,
for any $n \in {\mathbb N}$, and for an arbitrary $\varepsilon \in (0, \delta)$,
$$
\sup\limits_{|z|_p \le \delta - \varepsilon}\|y_n(z) - y(z)\| \to 0, \ n \to
\infty.
$$

Let $y \in {\mathfrak A}_{loc}(\mathfrak B)$. Its derivatives are defined as
$$
y^{(k)}(z) = \sum\limits_{n = 0}^{\infty} (n + 1) \dots (n + k) c_{n + k}z^n,
\quad k \in {\mathbb N}.
$$
It follows from (2) that
$$
r(y^{(k)}) \ge r(y).
$$
It is easily checked also that if $z \to 0$, then
$$
y(z) \to y(0) = c_0, \ \mbox{\rm and} \ \frac{y^{(k)}(z) - y^{(k)}(0)}{z}
\to y^{(k + 1)}(0) = c_{k + 1}(k + 1)!  \eqno (3)
$$
in the topology of $\mathfrak B$.

{\bf 3.} \ Let us consider the Cauchy problem
$$
\left\{
\begin{array}{rcl}
{\displaystyle \frac{dy(z)}{dz}} & = & Ay(z) \\
y(0) & = & y_0, \\
\end{array} \right.  \eqno (4)
$$
where $A$ is a closed linear operator on $\mathfrak B$. We say that a
vector-function $y(z)$ from ${\mathfrak A}_{loc}(\mathfrak B)$ is a solution of
problem (4) if $y(z) \in {\mathcal D}(A)$ for $z \in d(r(y)^-; \Omega)$ and satisfies
(4) in this disk.
\begin{teo}
In order that problem $(4)$ have a solution in ${\mathfrak A}_{loc}(\mathfrak B)$,
it is necesary and sufficient that $y_0 \in E(A)$; moreover
$\sigma(y_0; A)r(y) = {\displaystyle p^{-\frac{1}{p - 1}}}$. The problem $(4)$
is well-posed,
that is, the solution is unique, anf if a sequence of initial data
$y_{n,0} \in E(A)$ converges to $y_0$ in $E(A)$, then the sequence
of the corresponding solutions $y_n(z) \in {\mathfrak A}_{loc}(\mathfrak B)$ converges
to $y(z)$ in the space ${\mathfrak A}_{loc}(\mathfrak B)$.
\end{teo}
{\it Proof}. \ Suppose that
$$
y(z) = \sum\limits_{k = 0}^{\infty}c_kz^k, \quad c_k \in {\mathfrak B}, \quad
z \in d(r(y)^-; \Omega),
$$
is a solution of (4). Then $c_k \in {\mathcal D}(A), \ k \in {\mathbb N}_0$. Indeed,
$c_0 = y(0) = y_0 \in {\mathcal D}(A)$.
Since
$$
{\mathcal D}(A) \ni \frac{y(z) - c_0}{z} = \sum\limits_{k = 0}^{\infty}c_kz^{k - 1}
\to c_1,
$$
and
$$
\frac{A(y(z) - y_0)}{z} = \frac{y'(z) - y'(0)}{z}
\to 2c_2
$$
as $z \to 0$, we have because of closure of the operator $A$ that
$c_1 \in {\mathcal D}(A)$, and $Ac_1 = 2c_2$. Using (3), we get by induction that
$$
\forall k \in {\mathbb N} \quad c_k = \frac{Ac_{k - 1}}{k} \in {\mathcal D}(A),
$$
hence,
$$
\forall k \in {\mathbb N} \ y_0 \in {\mathcal D}(A^k), \ \mbox{\rm and} \ A^ky_0
= k!c_k.
$$
In view of (2),
$$
\frac{1}{r(y)} = \varlimsup_{n \to \infty}\sqrt[n]{\frac{\|A^ny_0\|}{|n!|_p}}.
$$
Taking into account the equality
$$
\lim\limits_{n \to \infty}\sqrt[n]{|n!|_p} = p^{-\frac{1}{p - 1}}
$$
(see [1]), we obtain
$$
\varlimsup_{n \to \infty}\sqrt[n]{\|A^ny_0\|} =
\frac{p^{-\frac{1}{p - 1}}}{r(y)},
$$
whence
$$
\forall n \in {\mathbb N}_0 \quad \|A^ny_0\| \le c{\alpha}^n,
$$
where $0 < c = \mbox{\rm const}, \ \alpha =
\frac{\displaystyle p^{-\frac{1}{p - 1}}}{\displaystyle r(y)}$.
So, $y_0 \in E(A)$, and $\sigma(y_0; A)r(y) = p^{-\frac{1}{p - 1}}$.

Conversely, let $y_0$ be an entire vector of exponential type for the operator
$A$ with $\sigma(y_0; A) = \sigma$. Then the series
$$
y(z) = \sum\limits_{k = 0}^{\infty}\frac{A^ky_0}{k!}z^k \eqno (5)
$$
is convergent in the disk $d(r^-; \Omega)$, where
$$
r = r(y) = \left(\varlimsup_{n \to \infty}\sqrt[n]{\frac{\|A^ny_0\|}{|n!|_p}}
\right)^{-1} = \frac{\lim\limits_{n \to \infty}\sqrt[n]{|n!|_p}}
{\varlimsup_{n \to \infty}\sqrt[n]{\|A^ny_0\|}} =
{\displaystyle \frac{p^{-\frac{1}{p - 1}}}{\sigma}}.
$$

We shall prove now that if $z \in d(r^-; \Omega)$, then $y(z) \in {\mathcal D}(A)$.
Really, since every component of series (5) belongs to ${\mathcal D}(A)$, the
sums $S_n(z) = \sum\limits_{k = 0}^{n}
{\displaystyle \frac{A^ky_0}{k!}z^k}$ belong to
${\mathcal D}(A)$, too. For $z \in d(r^-; \Omega)$, the sequence $S_n(z)$
converges to $y(z) \ (n \to \infty)$ in the topology of $\mathfrak B$. As
$\sigma (Ay_0; A) = \sigma (y_0; A) = \sigma$, the sequence $AS_n(z) =
\sum\limits_{k = 0}^{n}{\displaystyle \frac{A^{k + 1}y_0}{k!}z^k}, \ n \in
{\mathbb N}$, converges in $\mathfrak B \ (n \to \infty)$
in the same disk $d(r^-; \Omega)$. Since the operator $A$ is closed,
we have that $y(z) \in {\mathcal D}(A)$ when $z \in d(r^-; \Omega)$.

The formal differentiation of series (5) verifies that $y(z)$ satisfies
(4). Thus, the vector-function (5) is a solution of problem (4).

It remains to check the well-posedness of problem (4). Assume that
$y_{n, 0} \to y_0 \ (n \to \infty)$ in $E(A)$. This means that there exists
a number $\alpha > 0$ such that $y_{n, 0} \in E_{\alpha}(A)$ for sufficiently
large $n$, and $\|y_{n, 0} - y_0\|_{\alpha} \to 0$ as $n \to \infty$.
It follows from the above proof of sufficiency that
$$
r(y_n) \ge \frac{p^{-\frac{1}{p - 1}}}{\alpha}.
$$
So, we may take $\delta = {\displaystyle \frac{p^{-\frac{1}{p - 1}}}{\alpha}}$,
and to complete the proof, we need only show that for an arbitrary fixed
$\varepsilon \in (0, 1), \ \|y_n(z) - y(z)\| \to 0 \ (n \to \infty)$
uniformly in the disk $d((1 - \varepsilon)\delta^-; \Omega)$. We have
$$
\|y_n(z) - y(z)\| = \left\|\sum\limits_{k = 0}^{\infty}\frac{A^k(y_{n, 0} -
y_0)z^k}{k!}\right\| \le
\sum\limits_{k = 0}^{\infty} \frac{\|A^k(y_{n, 0} - y_0)\||z|_p^k}{|k!|_p} =
$$
$$
\sum\limits_{k = 0}^{\infty} \frac{\|A^k(y_{n, 0} - y_0)\|}{{\alpha}^k}\cdot
\frac{{\alpha}^k |z|_p^k}{|k!|_p} \le \sup\limits_{k}
\frac{\|A^k(y_{n, 0} - y_0)\|}{{\alpha}^k}\cdot \sum\limits_{k = 0}^{\infty}
\frac{{\alpha}^k |z|_p^k}{|k!|_p} \le \|y_{n, 0} - y_0\|_{\alpha}
\sum\limits_{k = 0}^{\infty} {\alpha}^k(1 - \varepsilon)^k
\frac{p^{-\frac{k}{p - 1}}}{{\alpha}^k|k!|_p}.
$$
Taking into account that
$$
\frac{1}{|k!|_p} \le {\displaystyle p^{\frac{k}{p - 1}}}
$$
(see [2]), we arrive at the inequality
$$
\|y_n(z) - y(z)\| \le {\varepsilon}^{-1}\|y_{n, 0} - y_0\|_{\alpha},
$$
and therefore problem (4) is well-posed. $\boxed{\phantom{\cdot}}$

It is seen from the proof of Theorem 1 that the series
$\sum\limits_{k = 0}^{\infty}{\displaystyle \frac{A^ky_0}{k!}z^k}$ is convergent for all
$z \in \Omega$ if and only if $y_0 \in \bigcap\limits_{\alpha > 0}E_{\alpha}(A)$.
\begin{cor}
If the operator $A$ is bounded, then for any $y_0 \in \mathfrak B$, the Cauchy
problem $(4)$ is well-posed in the class ${\mathfrak A}_{loc}(\mathfrak B)$.
\end{cor}
\begin{rem}
If $\mathfrak B$ is a Banach space over the field $\mathbb C$ of complex numbers,
then, as was shown in $[4]$, problem $(4)$ is well-posed in the class
of locally analytic vector-functions if and only if $y_0$ is an analytic vector
for the operator $A$, that is, $y_0 \in \bigcap\limits_{n \in {\mathbb N}_0}
{\mathcal D}(A^n)$, and
$$
\exists \alpha > 0 \quad \exists c > 0 \quad \forall k \in {\mathbb N}_0
\quad \|A^ky_0\| \le c{\alpha}^kk!.
$$
In this case, in order that problem $(4)$ be well-posed in the class
of entire vector-functions of exponential type, it is necessary and sufficient
that $y_0 \in E(A)$.
\end{rem}

{\bf 4.} \ In this section we show how the existence and uniqueness
theorem for the Cauchy problem for partial differential equations over a
non-archimedean field of characteristic zero (see [5]) may be
obtained from the above result.

Let ${\mathcal A}_{\rho}$ be the space of $\Omega$-valued functions $f(x)$ analytic
on the $n$-dimensional disk
$$
d({\rho}^+; {\Omega}^n) = \left\{x = (x_1, \dots x_n) \in {\Omega}^n:
|x|_p = \left(\sum\limits_{i = 1}^{n}|x_i|_p^2\right)^{1/2}
\le \rho\right\}.
$$
This means that
$$
f(x) = \sum\limits_{\alpha} f_{\alpha}x^{\alpha}, \quad f_{\alpha} \in \Omega,
\quad \lim\limits_{|\alpha| \to \infty}|f_{\alpha}|_p{\rho}^{|\alpha|} = 0,
$$
where $\alpha = ({\alpha}_1, \dots
{\alpha}_n), \ {\alpha}_i \in {\mathbb N}_0, \ |\alpha| = {\alpha}_1 + \dots + {\alpha}_n$.

The space ${\mathcal A}_{\rho}$ is a non-archimedean Banach space with respect
to the norm
$$
\|f\|_{\rho} = \sup\limits_{\alpha}|f_{\alpha}|_p {\rho}^{|\alpha|}.
$$
It is clear that the differential operators
$$
\frac{\partial f}{\partial x_j} = \sum\limits_{\alpha} {\alpha}_j f_{\alpha}
x_1^{{\alpha}_1} \dots x_j^{{\alpha}_j - 1} \dots x_n^{{\alpha}_n}, \ j =
1, \dots, n,
$$
are bounded in ${\mathcal A}_{\rho}$, and
$$
\left\|\frac{\partial f}{\partial x_j}\right\|_{\rho} =
\max\limits_{\alpha}|f_{\alpha}{\alpha}_j|_p{\rho}^{|\alpha| - 1} \le
\frac{1}{\rho}\max\limits_{\alpha}|f_{\alpha}|_p{\rho}^{|\alpha|} =
\frac{1}{\rho}\|f\|_{\rho}.
$$
The multiplication map
$$
G: f \mapsto fg,
$$
where $g \in {\mathcal A}_{\rho}$, is bounded in ${\mathcal A}_{\rho}$, too, and
$$
\|G\| = \|g\|_{\rho}.
$$
Indeed, let $f(x) = \sum\limits_{\alpha} f_{\alpha}x^{\alpha}, \
g(x) = \sum\limits_{\alpha} g_{\alpha}x^{\alpha}, \ \alpha = ({\alpha}_1,
\dots, {\alpha}_n)$. Then
$$
f(x)g(x) = \sum\limits_{\alpha} c_{\alpha}x^{\alpha},
$$
where
$$
c_{\alpha} = \sum\limits_{0 \le i \le \alpha} f_ig_{\alpha - i}
= \sum\limits_{i_1 = 0}^{{\alpha}_1} \dots \sum\limits_{i_n = 0}^{{\alpha}_n}
f_{i_1, \dots, i_n}g_{{\alpha}_1 - i_1, \dots, {\alpha}_n - i_n} \quad
(i = (i_1, \dots, i_n)).
$$
So,
$$
\|fg\|_{\rho} = \sup\limits_{\alpha}\max_{0 \le i \le \alpha} |f_i|_p
|g_{\alpha - i}|_p {\rho}^{|i|}{\rho}^{|\alpha| - |i|} \le \|f\|_{\rho}
\|g\|_{\rho}.
$$

We pass now to the Cauchy problem
$$
\left\{
\begin{array}{rcl}
{\displaystyle \frac{\partial u(t, x)}{\partial t}} & = &
\sum\limits_{|\beta| = 0}^na_{\beta}(x)D^{\beta}u(t, x) \\
u(0, x) & = & \varphi (x), \\
\end{array}
\right. \eqno (6)
$$
where $a_{\beta}(x) \in {\mathcal A}_{\rho}, \ \varphi (x) \in {\mathcal A}_{\rho}, \
D^{\beta} = {\displaystyle \frac{{\partial}^{|\beta|}}
{\partial x_1^{{\beta}_1} \dots \partial x_n^{{\beta}_n}}}$.

In the space ${\mathcal A}_{\rho}$ we define the operator $A$ as follows:
$$
f \mapsto Af = \sum\limits_{|\beta| = 0}^na_{\beta}D^{\beta}f.
$$
The relations
$$
\bigl\|\sum\limits_{|\beta| = 0}^na_{\beta}D^{\beta}f\bigr\|_{\rho} \le
\max\limits_{\beta}\|a_{\beta}D^{\beta}f\|_{\rho} \le
\max\limits_{\beta}\|a_{\beta}\|_{\rho}\|D^{\beta}f\|_{\rho} \le
\max\limits_{\beta}\{{\rho}^{-|\beta|}\|a_{\beta}\|_{\rho}\}\|f\|_{\rho}
$$
show that the operator $A$ is bounded, and
$$
\|A\| \le \max\limits_{\beta}\{{\rho}^{-|\beta|}\|a_{\beta}\|_{\rho}\}.
$$
It follows from Corollary 1 that problem (6) is well-posed in
${\mathfrak A}_{loc}({\mathcal A}_{\rho})$ in the disk
$$
\left\{{\displaystyle t \in \Omega: |t|_p <
\frac{p^{-\frac{1}{p - 1}}}
{\max\limits_{\beta}\{{\rho}^{-|\beta|}\|a_{\beta}\|_{\rho}\}}}\right\}.
$$

It should be noted that Theorem 1 is valid also in the case where $\mathfrak B$ is
a Banach space over an arbitrary non-archimedean field $K$ of characteristic
zero. But then the role of ${\displaystyle p^{-\frac{1}{p - 1}}}$ is played by a certain
constant $b = b(K)$ such that ${\displaystyle \frac{1}{|n!|_K} \le b^n}$ (see, for
instance, [6]).
\newpage

\vspace*{3mm}
Institute of Mathematics \\
National Academy of Sciences of Ukraine \\
3 Tereshchenkivs'ka \\
Kyiv 01601, Ukraine \\
E-mail: imath@horbach.kiev.ua

\end{document}